\documentclass[12pt]{amsart}
\usepackage{amsmath,amssymb,amsbsy,amsfonts,latexsym,amsopn,amstext,cite,amsxtra,euscript,amscd,bm,mathabx,,booktabs}
\usepackage{mathrsfs}
\usepackage{url}
\usepackage{amsbsy} 
\usepackage[colorlinks,linkcolor=blue,anchorcolor=blue,citecolor=blue,backref=page]{hyperref}
\usepackage{color}
\usepackage{graphics,epsfig}
\usepackage{graphicx}
\usepackage{float} 
\usepackage[english]{babel} 
\usepackage{mathtools}
\usepackage{todonotes}
\usepackage{url}

\usepackage[colorlinks,linkcolor=blue,anchorcolor=blue,citecolor=blue,backref=page]{hyperref}

\def\sL {\mathscr{L}}

\usepackage[norefs,nocites]{refcheck}

\hypersetup{breaklinks=true}

\usepackage[norefs,nocites]{refcheck}
\usepackage[english]{babel}
\begin{document}

\newtheorem{thm}{Theorem}
\newtheorem{lem}[thm]{Lemma}
\newtheorem{claim}[thm]{Claim}
\newtheorem{cor}[thm]{Corollary}
\newtheorem{prop}[thm]{Proposition} 
\newtheorem{definition}[thm]{Definition}
\newtheorem{rem}[thm]{Remark} 
\newtheorem{question}[thm]{Question}
\newtheorem{qn}[thm]{Question}
\newtheorem{conj}[thm]{Conjecture}
\newtheorem{prob}{Problem}
\newtheorem{assumption}[thm]{Assumption}

\newcommand{\GL}{\operatorname{GL}}
\newcommand{\SL}{\operatorname{SL}}
\newcommand{\lcm}{\operatorname{lcm}}
\newcommand{\ord}{\operatorname{ord}}
\newcommand{\Op}{\operatorname{Op}}
\newcommand{\Tr}{\operatorname{Tr}}
\newcommand{\Nm}{\operatorname{Nm}}

\numberwithin{equation}{section}
\numberwithin{thm}{section}
\numberwithin{table}{section}

\def\vol {{\mathrm{vol\,}}}
\def\squareforqed{\hbox{\rlap{$\sqcap$}$\sqcup$}}
\def\qed{\ifmmode\squareforqed\else{\unskip\nobreak\hfil
\penalty50\hskip1em\null\nobreak\hfil\squareforqed
\parfillskip=0pt\finalhyphendemerits=0\endgraf}\fi}

\def \balpha{\bm{\alpha}}
\def \bbeta{\bm{\beta}}
\def \bgamma{\bm{\gamma}}
\def \blambda{\bm{\lambda}}
\def \bchi{\bm{\chi}}
\def \bphi{\bm{\varphi}}
\def \bpsi{\bm{\psi}}
\def \bomega{\bm{\omega}}
\def \btheta{\bm{\vartheta}}
\def \ochi{\overline{\chi}}
\def \h{\widehat h}

\def\ovG {\overline{\Gamma}}

\def\uu {\mathbf{u}} 
\def\uv {\mathbf{v}} 
\def\uk {\mathbf{k}} 
\def\ut {\mathbf{t}}

\def\eps{\varepsilon}

\newcommand{\bfxi}{{\boldsymbol{\xi}}}
\newcommand{\bfrho}{{\boldsymbol{\rho}}}

\def\Kab{\sfK_\psi(a,b)}
\def\Kuv{\sfK_\psi(u,v)}
\def\SaUV{\cS_\psi(\balpha;\cU,\cV)}
\def\SaAV{\cS_\psi(\balpha;\cA,\cV)}

\def\SUV{\cS_\psi(\cU,\cV)}
\def\SAB{\cS_\psi(\cA,\cB)}

\def\Kmnp{\sfK_p(m,n)}

\def\KKap{\cH_p(a)}
\def\KKaq{\cH_q(a)}
\def\KKmnp{\cH_p(m,n)}
\def\KKmnq{\cH_q(m,n)}

\def\Klmnp{\sfK_p(\ell, m,n)}
\def\Klmnq{\sfK_q(\ell, m,n)}

\def \SALMNq {\cS_q(\balpha;\cL,\cI,\cJ)}
\def \SALMNp {\cS_p(\balpha;\cL,\cI,\cJ)}

\def \SACXMQX {\fS(\balpha,\bzeta, \bxi; M,Q,X)}

\def\SAMJp{\cS_p(\balpha;\cM,\cJ)}
\def\SAMJq{\cS_q(\balpha;\cM,\cJ)}
\def\SAqMJq{\cS_q(\balpha_q;\cM,\cJ)}
\def\SAJq{\cS_q(\balpha;\cJ)}
\def\SAqJq{\cS_q(\balpha_q;\cJ)}
\def\SAIJp{\cS_p(\balpha;\cI,\cJ)}
\def\SAIJq{\cS_q(\balpha;\cI,\cJ)}

\def\RIJp{\cR_p(\cI,\cJ)}
\def\RIJq{\cR_q(\cI,\cJ)}

\def\TWXJp{\cT_p(\bomega;\cX,\cJ)}
\def\TWXJq{\cT_q(\bomega;\cX,\cJ)}
\def\TWpXJp{\cT_p(\bomega_p;\cX,\cJ)}
\def\TWqXJq{\cT_q(\bomega_q;\cX,\cJ)}
\def\TWJq{\cT_q(\bomega;\cJ)}
\def\TWqJq{\cT_q(\bomega_q;\cJ)}

 \def \xbar{\overline x}
  \def \ybar{\overline y}

\def\cA{{\mathcal A}}
\def\cB{{\mathcal B}}
\def\cC{{\mathcal C}}
\def\cD{{\mathcal D}}
\def\cE{{\mathcal E}}
\def\cF{{\mathcal F}}
\def\cG{{\mathcal G}}
\def\cH{{\mathcal H}}
\def\cI{{\mathcal I}}
\def\cJ{{\mathcal J}}
\def\cK{{\mathcal K}}
\def\cL{{\mathcal L}}
\def\cM{{\mathcal M}}
\def\cN{{\mathcal N}}
\def\cO{{\mathcal O}}
\def\cP{{\mathcal P}}
\def\cQ{{\mathcal Q}}
\def\cR{{\mathcal R}}
\def\cS{{\mathcal S}}
\def\cT{{\mathcal T}}
\def\cU{{\mathcal U}}
\def\cV{{\mathcal V}}
\def\cW{{\mathcal W}}
\def\cX{{\mathcal X}}
\def\cY{{\mathcal Y}}
\def\cZ{{\mathcal Z}}
\def\Ker{{\mathrm{Ker}}}
\def\g{{\mathrm{gcd}}}
\def\rad{{\mathrm{rad}}}
\def\rank{{\mathrm{rank}}}
\def\Gal{{\mathrm{Gal}}}

\def\Psf {\mathsf P}
\def\Qsf {\mathsf Q}
\def\fsf {\mathsf f}
\def\ksf {\mathsf k}
\def\hsf {\mathsf h}
\def\isf{\mathsf O}

\def\Dsf {\mathsf D}

\def\NmQR{N(m;Q,R)}
\def\VmQR{\cV(m;Q,R)}

\def\Xm{\cX_m}

\def \A {{ \mathbb A}}
\def \B {{ \mathbb B}}
\def \C {{ \mathbb C}}
\def \N {{ \mathbb N}}
\def \F {{ \mathbb F}}
\def \G {{\mathbb G}}
\def \L {{\mathbb L}}
\def \K {{\mathbb K}}
\def \PP {{\mathbb P}}
\def \Q {{\mathbb Q}}
\def \R {{\mathbb R}}
\def \Z {{\mathbb Z}}
\def \fS{\mathfrak S}

\def\e{{\mathbf{\,e}}}
\def\ep{{\mathbf{\,e}}_p}
\def\eq{{\mathbf{\,e}}_q}
\def\er{{\mathbf{\,e}}_R}
\def\\{\cr}
\def\({\left(}
\def\){\right)}
\def\fl#1{\left\lfloor#1\right\rfloor}
\def\rf#1{\left\lceil#1\right\rceil}

\def\Tr{{\mathrm{Tr}}}
\def\Nm{{\mathrm{Nm}}}
\def\Im{{\mathrm{Im}}}

\def \oF {\overline \F}

\newcommand{\pfrac}[2]{{\left(\frac{#1}{#2}\right)}}

\def \Prob{{\mathrm {}}}
\def\e{\mathbf{e}}
\def\ep{{\mathbf{\,e}}_p}
\def\epp{{\mathbf{\,e}}_{p^2}}
\def\em{{\mathbf{\,e}}_m}

\def\Res{\mathrm{Res}}
\def\Orb{\mathrm{Orb}}

\def\vec#1{\mathbf{#1}}
\def \va{\vec{a}}
\def \vb{\vec{b}}
\def \vr{\vec{r}}
\def \vm{\vec{m}}
\def \vu{\vec{u}}
\def \vv{\vec{v}}
\def \vx{\vec{x}}
\def \vy{\vec{y}}
\def \vz{\vec{z}}

\newcommand{\gvec}{\boldsymbol{\gamma}}

\def\flp#1{{\left\langle#1\right\rangle}_p}
\def\T {\mathsf {T}}

\def\sfG {\mathsf {G}}
\def\sfK {\mathsf {K}}

\def\mand{\qquad\mbox{and}\qquad}

\title[Sparse systems of congruences]
{Counting solutions to sparse systems of congruences}

\author[S. Bhakta]{Subham Bhakta}
\address{School of Mathematics and Statistics, University of New South Wales, Sydney, NSW 2052, Australia.} 
\email{subham1729@gmail.com}

\author[I. E. Shparlinski]{Igor E. Shparlinski}
\address{School of Mathematics and Statistics, University of New South Wales, Sydney, NSW 2052, Australia.} 
\email{igor.shparlinski@unsw.edu.au}

\begin{abstract}
We obtain new bounds and asymptotic formulas for the number of solutions to sparse 
systems of $n$ diagonal congruences modulo an integer $q$ in $s$ variables. In a large range 
of parameters, our results improve and generalise those of A.~A.~Karatsuba  (1976). 
In particular, our lower bound on $s$,  which guarantees the solvability of such systems, 
depends only on $n$, rather than on specific degrees as in the work of A.~A.~Karatsuba  (1976).
\end{abstract}


\maketitle
\tableofcontents
\section{Introduction}

Let $\vr=(r_1,\ldots,r_n)\in \Z^n$ be a vector of pairwise distinct nonzero integers and let $q$
be a positive integer. We always assume that 
\begin{equation}
\label{eq:vec r}
|r_n|  \ge   \ldots  \ge  |r_1| \ge 1  \mand  \qquad r_1\cdots r_n \ne 0.
\end{equation}

Given a vector $\blambda= (\lambda_1, \ldots, \lambda_n) \in \Z_q^n$, we denote by $N_{\vr, s}(\blambda,q)$ the number of solutions to the symmetric 
system of congruences 
\begin{equation}
\label{eq:syst lambda}
\sum_{i=1}^s x_i^{r_j} \equiv\lambda_j\pmod q, \quad  1\leq j\leq n,\qquad x_1, \ldots, x_s,   \in \Z_q^*, 
\end{equation} 
in $s$  variables from the set $ \Z_q^*$ of   reduced residue classes 
in the residue ring $\Z_q$ modulo $q$.  

In the case when $q = p^k$ is a prime power,
systems of this type have been 
studied by Karatsuba~\cite{Kar},  see also the work of Wooley~\cite{Wool}.
In particular,  Karatsuba~\cite{Kar} shows that for $q =p^k$ and
positive integers  $r_n >  \ldots > r_1\ge 1$, the
system~\eqref{eq:syst lambda} has solutions, in variables from $\Z_q$ rather than $\Z_q^*$,
provided $s \ge 32 n r_n$.
Here, we improve the results of Karatsuba~\cite{Kar}, showing that the solvability condition on 
$s$ can be presented only in terms of the number of equations $n$ rather than of their degree, 
which is much stronger for sparse systems of equations. We also extend this result to arbitrary 
moduli $q$ and to negative exponents. Our approach builds on the ideas from~\cite{BhShpCat}.
 
In turn, our asymptotic results for $N_{\vr, s}(\blambda,q)$  depend on new bounds on 
$M_{\vr, s}(q)$,  which is the number of solutions to the symmetric 
system of congruences 
\begin{equation}
\label{eq:symm syst}
\begin{split}
\sum_{i=1}^s x_i^{r_j}\equiv\sum_{i=1}^s & y_i^{r_j}\pmod q, \qquad  1\leq j\leq n,\\
& x_1, y_1, \ldots, x_s, y_s \in \Z_q^*, 
\end{split}
\end{equation}
in $2s$  variables from the set $ \Z_q^*$ of   reduced residue classes 
in the residue ring $\Z_q$ modulo $q$.

To formulate our result, we  set $\gamma_{1,s} = 0$ and for integers $n\ge 2$ and $s\ge 1$, we define 
 \begin{equation}\label{eq:gamma_mh}
\gamma_{n,s} = 
 \min_{\substack{g_1,\ldots,g_{n-1}\in\Z\\
 g_1,\ldots,g_{n-1}\ge 0\\
                    g_1+\cdots+g_{n-1} =s-1}}\max\left\{0, 
\max_{1\leq \nu \leq n-1}
                   \sum_{j=\nu}^{n-1}\(1 -2 g_j/(j+1) \)\right\}, 
\end{equation}
 where the outer minimum is taken over all possible choices
 of non-negative integers $g_1,\ldots,g_{n-1}$ with $g_1+\cdots+g_{n-1} =s-1$.

Throughout the paper, the implied constants 
in  the relations
\[X = O(Y),  \qquad X \ll Y, \qquad Y \gg X
\] 
may depend only on $r_n$ and $s$. Furthermore, we write $V=U^{o(1)}$ for any quantity $V$, which for any $\delta > 0$ satisfies $|V| \le c(\delta) U^{\delta}$, with some constant $c(\delta)$ that may depend on $\delta$ (and both $r_n$ and $s$). In particular, we have $U^{o(1)} \cdot U^{o(1)} = U^{o(1)}$.

The letter $p$, always denotes a  prime number and $q$ denotes a positive integer.

\begin{thm}\label{thm: bound M} 
Let $\vr=(r_1,\ldots,r_n)$ be a fixed vector of $n \ge 1$ pairwise distinct nonzero 
integers, satisfying~\eqref{eq:vec r}.
  Then,  for any fixed $s$, we have 
 \[
M_{\vr, s}(q) 
\le    q^{2s-n+\gamma_{n,s} + o(1)},
\]
as $q\to \infty$, where $\gamma_{n,s}$ is given by~\eqref{eq:gamma_mh}.
\end{thm}

Next, using Theorem~\ref{thm: bound M}, we can obtain an asymptotic formula for $N_{\vr, s}(\blambda,q)$. 

  We denote by $\sL(q)$  the smallest prime divisor of $q$.
We also denote
 \begin{equation}\label{eq:Hn}
H_n=\sum_{k=1}^n k^{-1},
\end{equation}
and $\{x\}=x-\lfloor x\rfloor$ denotes the fractional part of $x$.
As usual, we use $\varphi(q)$ to denote the Euler function
of $q$.

\begin{thm}\label{thm: asymp N} 
Let $\vr=(r_1,\ldots,r_n)$ be a fixed vector of $n \ge 1$ pairwise distinct nonzero 
integers, satisfying~\eqref{eq:vec r}.
  Then for any fixed 
\[
s\ge \frac{n(n+3)}{2}+\lfloor 2H_n \rfloor, 
\]
uniformly over vectors $\blambda \in \Z_q^n$, we have 
 \[
N_{\vr, s}(\blambda,q) = \frac{\varphi(q)^s}{q^n} \(1 + 
O\(\sL(q)^{-\frac{1-\{2H_n\}}{n}+o(1)}\)\),
\] 
as $q\to \infty$, where $H_n$  is given by~\eqref{eq:Hn}.
\end{thm}

In particular, we see that $N_{\vr, s}(\blambda,q)> 0$ 
provided that $s$ satisfies the condition of Theorem~\ref{thm: asymp N} and $\sL(q)$ is large 
enough (in terms of only $r_n$ and $s$), improving and generalising the result of Karatsuba~\cite{Kar}
which which applies only for prime powers $q = p^k$ with $p> n$,  assumes the positivity $r_n > \ldots > r_1\ge 1$
and, more importantly, requires $s$ to be large enough in terms of $r_n$ rather than of $n$.

On the other hand, we also note that~\eqref{eq:sn upper lower} shows that our approach does not allow any $s\le n(n+3)/2+1$.

In fact, if we restrict to the case of prime powers $q = p^k$ as in Karatsuba~\cite{Kar}, we get the existence of solutions to~\eqref{eq:syst lambda},  that is, $N_{\vr, s}(\blambda,q)> 0$, for even smaller $s$, provided $p$ is large enough (in terms of only $r_n$ and $s$). More precisely, we prove the following.

\begin{thm}\label{thm: asymp pk} 
Let $q = p^k$ be a prime power and $\vr=(r_1,\ldots,r_n)$ be a fixed vector of $n \ge 1$ pairwise 
distinct nonzero integers, satisfying~\eqref{eq:vec r}. Then for any fixed 
\[
s\ge \frac{n(n+1)}{2}+\lfloor 2H_n \rfloor, 
\]
uniformly over vectors $\blambda \in \Z_q^n$, we have 
 \[
N_{\vr, s}(\blambda,q) = \frac{\varphi(q)^s}{q^n} \(1 + 
O\(p^{-\frac{1-\{2H_n\}}{n}+o(1)}\)\),
\] 
as $p\to \infty$, where $H_n$  is given by~\eqref{eq:Hn}.
\end{thm}

We emphasise that the implied constant in the error term of Theorem~\ref{thm: asymp pk}
does not depend on $k$, and thus there exist some $p_0(r_n,s)$ depending only 
on $r_n$ and $s$ that for $p \ge p_0(r_n,s)$, for all $k =1,2, \ldots$, we have 
\[
N_{\vr, s}(\blambda,q) \ge  \frac{\varphi(q)^s}{2q^n} > 0. 
\]

Finally, we remark that our results extend without any changes, albeit some obvious trypographical adjustments, 
to the systems of congruneces similar to~\eqref{eq:syst lambda}, but with arbitrary coefficients relatively prime to $q$, that is, 
with $x_i^{r_j}$ replaces with $a_{i,j} x_i^{r_j}$, $a_{i,j} \in \Z$, $\gcd\(a_{i,j}, q\)=1$, $i =1, \ldots, s$, $j=1, \ldots, n$. A similar comment also applies to~\eqref{eq:symm syst}.

\section{Asymptotic behavior of $\gamma_{n,s}$}
\label{sec:num asymp gamma}

Let us consider
\[
m=\min\left\{n-1,\left\lfloor(4s+3)^{1/2}\right\rfloor-2\right\},
\]
and set 
$$M=s-\left\lfloor (m+2)^2/4\right\rfloor \ge 0.$$
Moreover, we denote $g_1=M+1$ and 
\[g_j=
\begin{cases}
\left\lceil (j+1)/2\right\rceil, & \mathrm{if }\quad 2\le j\le m,\\
0, & \mathrm{if }\quad m<j\le n-1.
\end{cases}
\]

Then we have
\[
g_1+\cdots+g_{n-1} =s-1,
\] 
and also 
\[
\max\left\{0,\max_{1\le \nu\le n-1}\sum_{j=\nu}^{n-1}\left(1-2g_j/(j+1)\right)\right\}
=n-1-m.\]

Consequently, we have a general upper bound
\[
\gamma_{n,s}\le \max\left\{0,n+1-\left\lfloor(4s+3)^{1/2}\right\rfloor\right\},
\]
and in particular,
\[
\gamma_{n,s}=0,\quad \mathrm{for}\quad s\ge \left\lfloor (n+1)^2/4\right\rfloor,
\]
which provides a better understanding of the savings in Theorem~\ref{thm: bound M}.

For any integer $n\ge 1$, let us now denote
\begin{equation}\label{eq:An}
A_n=\min_{v=1,2,\ldots}\(2v+n\gamma_{n,v}\),
\end{equation}
and set
\begin{equation}\label{eq:s to A}
s(n)=n+\lfloor A_n\rfloor+1.
\end{equation}

Now we take a positive integer $u$ (of course, depending on $n$) attaining the minimum in the definition of $A_n$. Then we have 
\begin{equation}\label{eq:u of An}
A_n=2u+n\gamma_{n,u}.
\end{equation}
In particular we have $u< s(n)/2$, as  
\[
2u\le 2u+n\gamma_{n,u}=A_n<s(n)-n<s(n).
\]

On the other hand, since we have
\[
\gamma_{n,u}=\min_{\substack{g\ge 0, h\ge 1\\ h+g=u}}\max\left\{0,1+\gamma_{n-1,h}-2g/n\right\},
\]
consequently we derive
\begin{equation}\label{eq:An vs gamma n-1}
\begin{split}
A_n &= \min_{g\ge 0, h\ge 1}\max\left\{2(h+g),n+2h+n\gamma_{n-1,h}\right\}\\
&= n+\min_{h\ge 1}\(2h+n\gamma_{n-1,h}\).
\end{split}
\end{equation}

For any any integer $n\ge 1$, let us write
\[A_n=2q_n+\rho_n,\quad\mathrm{where}\quad q_n=\left\lfloor A_n/2\right\rfloor\quad \mathrm{and}\quad 0\le \rho_n <2.
\]

\begin{lem}\label{lem:An recurrence}
The sequence $\(A_n\)_{n \ge 1}$ satisfy the following recurrence relation
\[
A_{n+1}=n+1+A_n+\min\left\{\rho_n/n,2-\rho_n\right\},\qquad n=1,2,\ldots.
\]
\end{lem}

\begin{proof}
We first  show that
\begin{equation}\label{eq:gamma-endpoint-values}
\gamma_{n,q_n}=\rho_n/n\quad \mathrm{and}\quad \gamma_{n,q_n+1}=0.
\end{equation}

Again, let $1\le u< s(n)/2$ be an integer  $v$,  attaining the minimum in the definition of $A_n$. So we can write
\[
A_n=2u+n\gamma_{n,u}.
\]
Since $\gamma_{n,u}\ge0$, we also have $u\le A_n/2$, and in particular, $u\le q_n$.

On the other hand, for any integer $t\ge 0$, trivially we have
\[\gamma_{n,u+t}\le
\max\left\{0,\gamma_{n,u}-\frac{2t}{n}\right\}=
\max\left\{0,\frac{A_n-2(u+t)}{n}\right\}.\] 


In particular, we have 
\[
\gamma_{n,q_n} \le\frac{A_n-2q_n}{n}=\rho_n/n.
\] 
On the other hand, the equation~\eqref{eq:An} implies that $A_n\le 2q_n+n\gamma_{n,q_n}$, and hence $\gamma_{n,q_n}\ge \rho_n/n$, indeed this shows that $\gamma_{n,q_n}=\rho_n/n$. 

Also taking $t=q_n+1-u$, it follows that
\[
0\le\gamma_{n,q_n+1}
\le
\max\left\{0,(\rho_n-2)/n\right\}
=0,
\]
and this completes the proof of~\eqref{eq:gamma-endpoint-values}.

On the other hand, for any $1\le h\le q_n$, the definition of $A_n$ implies
\[
\gamma_{n,h}\ge\frac{A_n-2h}{n},
\]
and hence
\[
2h+(n+1)\gamma_{n,h} \ge A_n+\rho_n/n.
\]
Moreover by~\eqref{eq:gamma-endpoint-values}, the equality is attained at $h=q_n$, and hence
\begin{equation}\label{eq:h le q}
\min_{1\le h\le q_n} \left(2h+(n+1)\gamma_{n,h}\right)
=A_n+\rho_n/n.
\end{equation}

For $h\ge q_n +1$, the non-negativity of $\gamma_{n,h}$ implies
\[
2h+(n+1)\gamma_{n,h}\ge 2h \ge 2(q_n+1)=A_n+2-\rho_n.
\]
Again by~\eqref{eq:gamma-endpoint-values}, the equality is attained at $h=q_n+1$, and hence
\begin{equation}\label{eq:h ge q+1}
\min_{h\ge q_n+1} \left(2h+(n+1)\gamma_{n,h}\right)=A_n+2-\rho_n.
\end{equation}

The proof finally concludes by combining~\eqref{eq:An vs gamma n-1}, \eqref{eq:h le q} and~\eqref{eq:h ge q+1}.
\end{proof}

Since $0\le \rho_n<2$, it is easy to see that
\[
0\le \min\left\{\rho_n/n,2-\rho_n\right\}\le 2/(n+1),\quad n=1,2,\ldots.
\] 
Moreover note that $A_1=2$, as $\gamma_{1,v}=0$ for any $v\ge 1$. In particular, Lemma~\ref{lem:An recurrence} implies
\begin{equation}\label{eq:An upper lower bound}
\frac{n(n+1)}{2}+1\le A_n \le \frac{n(n+1)}{2}-1+2H_n.
\end{equation}
Then it follows from~\eqref{eq:s to A} that
\begin{equation}\label{eq:sn upper lower}
\frac{n(n+3)}{2}+2\le s(n)\le \frac{n(n+3)}{2}+ \lfloor 2H_n\rfloor.
\end{equation}

Therefore~\eqref{eq:u of An} and~\eqref{eq:An upper lower bound} show that, for any integer $s\ge s_0(n)$
where 
\begin{equation}\label{eq:min s needed}
 s_0(n) = \frac{n(n+3)}{2}+\lfloor 2H_n \rfloor,
\end{equation}
we have
\begin{equation}\label{eq:exp of LQ}
2u/n + \gamma_{n,u} - s/n +1 \le  \(A_n-s_0(n)\)/n+1\le -\frac{1-\{2H_n\}}{n}.
\end{equation}

\section{Exponential sums with sparse polynomials}

Let us denote $\eq(z) = \exp(2\pi i z/q)$. For vectors 
 $\vr=(r_1,\ldots,r_n)\in \Z^n$ and $\va= ( a_1, \ldots, a_n) \in \Z_q^n$, we consider the sums
 \begin{equation}\label{eq:Sqar}
S_q(\va, \vr) =   \sum_{x\in  \Z_q^*} 
 \eq\(a_1x^{r_1} + \ldots + a_{n}x^{r_n}\).
\end{equation}

The following bound is given in~\cite[Theorem~1]{Shp}
(see also~\cite{Bourg, BhShp} for other bounds on exponential sums with sparse polynomials). 

\begin{lem}\label{lem: exp sum} 
Let $\vr=(r_1,\ldots,r_n)$ be a fixed vector of $n \ge 1$ pairwise distinct nonzero integers, satisfying~\eqref{eq:vec r}. Then, uniformly over vectors $\va \in \Z_q^n$ with 
$\gcd(\va, q) =  1$, we have 
 \[
\left| S_q(\va, \vr) \right| \le q^{1-1/n + o(1)},
\]
as $q\to \infty$. 
\end{lem}

\section{Proof of Theorem~\ref{thm: bound M}}

We prove the result by induction on $s$

For $s=1$ we recall that for a fixed integer $r \ne 0$,  by~\cite[Lemma~5]{FHS} the number of 
solutions to the congruence 
\[
x^r \equiv y^r \pmod q
\]
is at most $q^{1+o(1)}$.

We now assume that the result holds for all fixed vectors 
\[
\vr=(r_1,\ldots,r_h) \in \Z^h, \qquad r_i \ne r_j \ne 0, \quad i,j =1, \ldots h, \ i \ne j,
\]
of dimension $h \le n-1$ and establish it for vectors $\vr$ of dimension $n$. 

By orthogonality of additive characters, let us first write
  \begin{equation}\label{eq:R = S(a)}
   M_{\vr, s}(q)  
   =  \frac{1}{q^{n}}\sum_{\va \in \Z_{q}^{n}}\left|S_q(\va, \vr)  \right|^{2s},
\end{equation}
where $S_q(\va, \vr)$ be as in~\eqref{eq:Sqar}.

First, for each divisor $d$ of $q$, we collect together the vectors with 
\[
\gcd(\va, q) =  q/d
\]
and then for each such $d$, we have
\begin{equation}\label{eq: Sq Sd}
S_q(\va, \vr)  = \frac{\varphi(q)}{\varphi(d)} S_d(\va (q/d)^{-1}, \vr).
\end{equation}

Simply noting that $\varphi(q)\le  \varphi(d) q/d$, and writing $\va (q/d)^{-1} = \vb$ for $\vb = (b_1, \ldots, b_n) \in  (\Z_{d})^{n}$ with $\gcd(\vb, d) = \gcd(b_1, \ldots, b_{n}, d) =  1$, we  derive from~\eqref{eq:R = S(a)} and~\eqref{eq: Sq Sd} that  
\begin{equation}\label{eq: M va Sd}
    M_{\vr, s}(q)   \le q^{2s-n}   +  \frac{1}{q^n} \sum_{\substack{d\mid q\\ 1< d\le q}} 
     (q/d)^{2s}    \sum_{\substack{\vb \in (\Z_{d})^{n}\\ \gcd(\vb, d)=1}}\left| S_d(\vb, \vr) \right|^{2s}. 
\end{equation}

Next, for any integers $g\ge 0$ and $h\ge 1$ with $g+h = s$, let us write
\[
\left|  S_d(\vb, \vr) \right|^{2s} = \left|  S_d(\vb, \vr) \right|^{2g}\left|  S_d(\vb, \vr) \right|^{2h}.
\]
Then by Lemma~\ref{lem: exp sum}, we have
\begin{equation}\label{eq:Mrs to Mrh}
\begin{split}
   & M_{\vr, s}(q)  \\
    &\quad \le q^{2s-n}   +  \frac{1}{q^n} \sum_{\substack{d\mid q\\ 1< d\le q}} 
     (q/d)^{2s}  d^{2g(1-1/n+o(1))}  \sum_{\substack{\vb \in (\Z_{d})^{n}\\ \gcd(\vb, d)=1}}\left| S_d(\vb, \vr) \right|^{2h}\\
     &\quad \le  q^{2s-n}+ q^{2s+o(1)}   \sum_{\substack{d\mid q\\ 1< d \le q}}   d^{2g(1-1/n)-2s}  M_{\vr, h}(d),
\end{split}
\end{equation} 
where in the second inequality in~\eqref{eq:Mrs to Mrh}, we use the orthogonality relation~\eqref{eq:R = S(a)}, with the crude inequality $d\le q$ and by trivially extending the summation over $\vb$ with $\gcd(\vb, d)=1$ to all vectors $\vb \in (\Z_{d})^{n}$.

Now we denote $\vr^* = \(r_1, \ldots, r_{n-1}\)$. Discarding the $n$th congruence, trivially we have
\begin{equation}\label{eq:r to r*}
 M_{\vr, h}(d) \le  M_{\vr^*, h}(d)\le d^{2h-n+1+\gamma_{n-1,h} + o(1)},
\end{equation}
for any divisor $d$ of $q$.

Applying the inductive assumption, combining~\eqref{eq:Mrs to Mrh} and~\eqref{eq:r to r*}, we finally obtain
\begin{align*}
   M_{\vr, s}(q)  & \le   q^{2s-n}+ q^{2s+o(1)}   \sum_{\substack{d\mid q\\ 1< d \le q}}   d^{2g(1-1/n)-2s}  d^{2h-n+1+\gamma_{n-1,h}}\\
& \le q^{2s-n}+ q^{2s-n+\gamma_{n-1,h}+1-2g/n+o(1)},
\end{align*} 
where in the last inequality above, we use the standard estimate that the number of distinct positive integer factors of $q$, is at most $q^{o(1)}$; see for instance~\cite[Theorem~317]{HaWr}.

 

We now choose $g$ and $h$  to realise the minimum in 
\[
 \gamma_{n,s}  = \min_{\substack{g \ge 0, \ h \ge 1\\ g+h = s}} \max\{0, \gamma_{n-1,h} +1- 2g/n\},
\]
and this completes the proof. 

\section{Proofs of Theorems~\ref{thm: asymp N} and~\ref{thm: asymp pk}}

As in~\eqref{eq:Mrs to Mrh}, by orthogonality of additive characters, we write
\[
N_{\vr, s}(\blambda,q) =  
\frac{1}{q^{n}}\sum_{\va \in \Z_{q}^{n}} \(S_q(\va, \vr)  \)^{s} \eq\(-\langle \va, \blambda \rangle\).
\]
Seprating the contribution from the zero vector $\va = \mathbf 0 \in (\Z_{d})^{n}$ and using~\eqref{eq: Sq Sd} that  
we obtain the following version of~\eqref{eq: M va Sd} (where we keep the original ratio $\frac{\varphi(q)}{\varphi(d)}$
without replacing it with $q/d$):
\begin{equation}\label{eq: N va Sd}
\left| N_{\vr, s}(\blambda,q) -   \varphi(q)^{s}/q^{n}  \right| \le  \frac{1}{q^n} \sum_{\substack{d\mid q\\ 1< d\le q}} 
     \(\frac{\varphi(q)}{\varphi(d)}\)^{s}    \sum_{\substack{\vb \in (\Z_{d})^{n}\\ \gcd(\vb, d)=1}}\left| S_d(\vb, \vr) \right|^{s}. 
\end{equation}
We now choose an arbitrary integer $u <s/2$ and set $v = s-2u$. Using Lemma~\ref{lem: exp sum}  
and then Theorem~\ref{thm: bound M} we obtain 
\begin{align*}
\sum_{\substack{\vb \in (\Z_{d})^{n}\\ \gcd(\vb, d)=1}}\left| S_d(\vb, \vr) \right|^{s}& \le 
d^{v(1-1/n) + o(1)}  \sum_{\substack{\vb \in (\Z_{d})^{n}\\ \gcd(\vb, d)=1}}\left| S_d(\vb, \vr) \right|^{2u}\\
& \ll d^{v(1-1/n) + o(1)}  \sum_{\vb \in (\Z_{d})^{n}}\left| S_d(\vb, \vr) \right|^{2u} \\
& =  d^{v(1-1/n) + n + o(1)} M_{\vr, 2u}(d)  \\
& =  d^{v(1-1/n) + 2u+\gamma_{n,u} + o(1)} = d^{s - v/n + \gamma_{n,u} + o(1)} . 
\end{align*}

Substituting in~\eqref{eq: N va Sd} we obtain
\begin{equation}\label{eq:error for q}
\begin{split}
\left| N_{\vr, s}(\blambda,q) -   \varphi(q)^{s}/q^{n}  \right| & \le   \varphi(q)^{s} q^{- n}  \sum_{\substack{d\mid q\\ 1< d\le q}} (d/\varphi(d))^s
d^{- v/n + \gamma_{n,u}  +o(1)} \\
& = 
 \varphi(q)^{s} q^{- n}   \sum_{\substack{d\mid q\\ 1< d\le q}}  
d^{- v/n + \gamma_{n,u}  +o(1)}\\
& = 
 \varphi(q)^{s} q^{- n}   \sum_{\substack{d\mid q\\ 1< d\le q}}  
d^{2u/n + \gamma_{n,u} - s/n +o(1)}. 
\end{split}
\end{equation}

For an arbitrary integer $q$, we now derive 
\begin{equation}\label{eq:error for comp q}
\left| N_{\vr, s}(\blambda,q) -   \varphi(q)^{s}/q^{n}  \right|   \le \varphi(q)^{s} q^{- n}  \sL(q) ^{2u/n + \gamma_{n,u} - s/n +1 +o(1)}, 
\end{equation}
provided $2u/n + \gamma_{n,u} - s/n +1 < 0$.

The proof is now complete from~\eqref{eq:exp of LQ}, as the value of $u$ appearing in~\eqref{eq:u of An} satisfies
\[
2u< s(n) \le s_0(n) \le s,
\]
where $s(n)$ and $s_0(n)$ are defined by~\eqref{eq:s to A} and~\eqref{eq:min s needed}, respectively, and thus is admissible, which concludes the proof of Theorem~\ref{thm: asymp N}.

Let us now denote
\[\widetilde s(n)=\lfloor A_n\rfloor+1.\]
Again, for the same $u$ appearing in~\eqref{eq:u of An}, we have $$2u\le A_n<\widetilde s(n)\le s_0(n)-n.$$ 
In particular for any $s\ge s_0(n)-n$, it follows from~\eqref{eq:exp of LQ} that
\[2u/n + \gamma_{n,u} - s/n=\(A_n-s_0(n)\)/n+1\le -\frac{1-\{2H_n\}}{n}.\]

We note that for $q = p^k$, instead of~\eqref{eq:error for comp q} we derive from~\eqref{eq:error for q} 
a stronger bound 
\begin{align*}
\left| N_{\vr, s}(\blambda,q) -   \varphi(q)^{s}/q^{n}  \right| & \le   \varphi(q)^{s} q^{-n}   \sum_{1\le i\le k}  
p^{k\(2u/n + \gamma_{n,u} - s/n +o(1)\)}\\
&  \le \varphi(q)^{s} q^{-n} p^{k\(2u/n + \gamma_{n,u} - s/n +o(1)\)},
\end{align*}
which now implies Theorem~\ref{thm: asymp pk}.

\section*{Acknowledgement} 
During the preparation of this work, the authors were supported in part by the Australian Research Council Grant DP230100530. I.S. would also like to thank the Max Planck Institute for Mathematics, Bonn, for its hospitality, stimulating environment, and generous support. 

ChatGPT-6 has been used to derive the asymptotic results in Section~\ref{sec:num asymp gamma}; 
the other arguments and the exposition are due to the authors.

\end{document}